\documentclass[12pt]{article} 

 \usepackage{latexsym,amsfonts,amssymb}
      \usepackage{amsmath, esint}
      \newtheorem{theorem}{Theorem}[section]

      \newtheorem{remark}{Remark}[section]
      
      \def\nn{\nonumber}
      \def\rf#1{\mbox{$(\ref{#1})$}}

      \def\be{\begin{equation}} 
      \def\ee{\end{equation}} 
      \def\beqn{\begin{eqnarray}} 
      \def\eeqn{\end{eqnarray}} 
      \def\beq{\begin{eqnarray*}} 
      \def\eeq{\end{eqnarray*}}
      \def\proof{{\bf Proof:}\ }
      \def\ga{\gamma} 
      \def\ep{\varepsilon} 
      \def\la{\lambda} 
      \def\ra{\rightarrow} 
      
      \def\ack#1{{\bf Acknowledgements.}\mbox{}#1}

\begin{document}

\title{\bf Hierarchical Dirichlet Process and Relative Entropy\thanks{Supported
    by the Natural Sciences and Engineering Research Council of Canada.}}
    
     \author{Shui Feng\\Department of Mathematics and Statistics\\McMaster
      University\\Hamilton, Ontario
      Canada L8S 4K1\\
      shuifeng@mcmaster.ca}
      \date{\empty}
      \maketitle

\begin{abstract}
The Hierarchical Dirichlet process is a discrete random measure serving as an important prior in Bayesian non-parametrics. It is motivated with the study of groups of clustered data.  Each group is modelled through a level two Dirichlet process and all groups share the same base distribution which itself is a drawn from a level one Dirichlet process. It has two concentration parameters with one at each level. The main results of the paper are the law of large numbers and large deviations for the hierarchical Dirichlet process and its mass when both concentration parameters converge to infinity.  The large deviation rate functions are identified explicitly. The rate function for the hierarchical Dirichlet process consists of two terms corresponding to the relative entropies at each level. It is less than the rate function for the Dirichlet process, which reflects the fact that the number of clusters under the hierarchical Dirichlet process has a slower growth rate than under the Dirichlet process.

      \vspace*{.125in} \noindent {\bf Key words:} Dirichlet process, Poisson-Dirichlet
      distribution, hierarchical Dirichlet process, stick-breaking, subordinator, Bayesian non-parametrics.
      \vspace*{.125in}

      \noindent {\bf MSC2020-Mathematics Subject Classifications:}
      Primary: 60G57; Secondary: 62F15.

\end{abstract}

\newcommand{\ABS}[1]{\left(#1\right)} 
\newcommand{\veps}{\varepsilon} 




    \section{Introduction}
      \setcounter{equation}{0}

The {\it Dirichlet process} introduced in \cite{Fer73} is a random discrete probability  that serves as a fundamental prior in Bayesian nonparametrics. It has two parameters, a concentration parameter $\alpha$ and a base probability distribution $\nu_0$.   The base distribution is the prior mean, and the concentration parameter is inversely proportional to the prior variance.  Due to its discrete nature the Dirichlet process is an effective prior in inferences for clustering data.  The Dirichlet process with infinite concentration parameter is simply $\nu_0$,  which corresponds to the classical parametric setting.

In \cite{TJBB06},  the authors introduced the {\it Hierarchical Dirichlet process} (henceforth HDP). It is motivated for the study of groups of clustered data where each group is modelled through a Dirichlet process and all groups share the same base distribution which itself is a drawn from another Dirichlet process. The special sharing mechanism makes the model an ideal prior for data with more concentrated clusters (\cite{BGQ21}, \cite{EDFAT19}). The HDP has three parameters, the level one and level two concentration parameters, and the base distribution.

In addition to Bayesian statistics, hierarchical models are also natural and fundamental in many other  areas.  One notable situation is in statistical mechanics where it is used to describe interactions at different scales in a physical system (\cite{CoEc78}, \cite{DGW04}, \cite{DaiFP21}). Various studies on asymptotic behaviours have led to deep understandings of the scale specific interactions and interactions between different scales (\cite{DGW04}, \cite{DaGrVa95}).  In the same spirit, we investigate the asymptotic behaviours of the HBP. The limiting procedures involve large concentration parameters at both levels. 
 
The asymptotic behaviour of the Dirichlet process has been studied extensively in the last thirty years (\cite{Feng10}, \cite{GhVa17} and references therein). The hierarchical structure in HBP presents new challenges. The mathematical framework is similar to random systems in a random environment. Our results will reveal explicitly the level-specific and cross level impact of the corresponding randomness.

 The basic setting will be presented in Section 2. This will include necessary notation, the relative entropy, the Dirichlet process,  the HDP, and the related asymptotic results. The main results will be discussed in Section 3 including the law of large numbers and large deviations.
             
\section{Preliminaries}
\setcounter{equation}{0}

  Let $(\Omega,{\cal F}, P)$ be a probability space, $E=[0,1]$, and ${\cal E}$ be the $\sigma$-algebra of Borel subsets of $E$. Let $C(E)$ and $B(E)$ denote the spaces of  continuously functions and bounded measurable functions on $E$, respectively. Let $M_1(E)$ denote the  space of probability measures equipped with the weak topology and the metric
  \[
  \rho(\mu,\nu)=\sum_{i=1}^{\infty}\frac{|\langle \mu-\nu, f_i\rangle|\wedge 1}{2^i}
  \]
  where $\{f_i: i\geq 1\}$ is a dense subset of $C(E)$.

For any $\nu_0$ in $M_1(E)$ and $\alpha >0$, let $\xi_1, \xi_2, \ldots$ be i.i.d. with common distribution $\nu_0$, and independently $U_1, U_2, \ldots$ be i.i.d. with Beta$(1, \alpha)$ distribution.  Set
\be\label{gem}
V_1=U_1,\ V_n =(1-U_1)\cdots(1-U_{n-1})U_n, \ n \geq 2
\ee
and
\[
\mathbb{V}_\alpha=(V_1,V_2, \ldots).
\]

The Dirichlet process with concentration parameter $\alpha$ and base distribution $\nu_0$ is  given by 
\be\label{dp}
\Xi_{\alpha,\nu_0}=\sum_{i=1}^{\infty}V_i \delta_{\xi_i}.
\ee 

The Hierarchical Dirichlet process introduced in \cite{TJBB06} is a non-parametric model for the study of groups of data. The prior for each group is a (level two) Dirichlet process and all groups share the same base measure which itself is a draw from another (level one) Dirichlet process (hence the hierarchical structure). Given the common base measure, the Dirichlet processes for different groups are independent and identically distributed.  Since all Dirichlet processes for different groups share the same types as the level one Dirichlet process, a stick breaking representation is also obtained in \cite{TJBB06}.   

More specifically,  for $\beta>0$ and any $n\geq 1$, let $W_n$ be a Beta$(\beta V_n, \beta(1-\sum_{k=1}^n V_k))$ random variable. The random variables $W_1, W_2, \ldots $ are conditionally independent given $\mathbb{V}_\alpha$.  Define
\be\label{gem-r}
Z_1 =W_1,\ Z_n =(1-W_1)\cdots(1-W_{n-1})W_n, \ n\geq 2
\ee
and 
\[
\mathbb{Z}_{\alpha,\beta}=(Z_1, Z_2, \ldots).
\]

The HDP with level two concentration parameter $\beta$, level one concentration parameter $\alpha$, and base distribution $\nu_0$ is the random measure 

\be\label{hdp}
\Xi_{\alpha,\beta, \nu_0} \stackrel{d}{=} \Xi_{\beta, \Xi_{\alpha, \nu_0}}\stackrel{d}{=}\sum_{i=1}^{\infty} Z_i\delta_{\xi_i},
\ee
where $\stackrel{d}{=} $ denotes equality in distribution. The first equality is by definition (Dirichlet with Dirichlet base) and the second equality is the stick breaking representation.
  
Let $Q_{\alpha}$ and $Q_{\alpha,\beta}$ denote the respective distributions of $\mathbb{V}_\alpha$ and $\mathbb{Z}_{\alpha,\beta}$. The distributions of $\Xi_{\alpha,\nu_0}$ and $\Xi_{\alpha,\beta, \nu_0}$ will be denoted by $\Pi_{\alpha,\nu_0}$ and $\Pi_{\alpha,\beta,\nu_0}$ respectively.  When the concentration parameters tend to infinity, the prior will concentrate on the base distribution $\nu_0$ and the nonparametric model becomes parametric model. To understand the microscopic transition between these two types of models, it is natural to investigate the asymptotic behaviour associated with these limiting procedures. This includes the law of large numbers and the large deviations. 

The family of probability measures $\{P_{\la}:\la>0\}$ on a Polish space $S$ satisfies a {\it large deviation principle} with {\it speed} $a(\la)$ and  {\it good rate function} $I(\cdot)$ as $\la$ tends to infinity if 
\[
\{s\in S: I(s)\leq c\} \ \mbox{is compact for all finite }\ c\geq 0
\] 
and 
\[
-\inf_{s\in G^{\circ}}I(s) \leq \liminf_{\la \ra \infty}\frac{1}{a(\la)}\log P_\la(G) \leq \limsup_{\la \ra \infty}\frac{1}{a(\la)}\log P_\la(G)\leq -\inf_{s\in \bar{G}}I(s)\]
where $a(\la)$ converges to infinity, $G^{\circ}$ and $\bar{G}$ denote the interior and closure of the measurable set $G\subset S$ respectively.

The large deviations obtained in this paper will be for the families $\{Q_{\alpha,\beta}: \alpha>0, \beta>0\}$ and $\{\Pi_{\alpha,\beta}: \alpha>0, \beta>0\}$. It turns out that the large deviation results depend on the relative growth magnitude of the concentration parameters.  To capture and to compare the impact of both levels of randomness on large deviations we will focus on the limiting procedure:

\be\label{assump}
\alpha \ra \infty,\ \beta \ra \infty, \ \frac{\alpha}{\beta} \ra c \in (0,\infty).
 \ee

A main quantity for our large deviation results is the relative entropy or the Kullback-Leibler divergence.  More specifically, for any two probabilities $\mu$ and $\nu$ in $M_1(E)$ the relative entropy of $\mu$ with respect to $\nu$ is defined and denoted by
\be\label{re-entropy}
H(\mu|\nu)=\left\{\begin{array}{ll}
  \int_E \log \frac{d\,\mu}{d\,\nu}\ \mu(d\,x)&  \mu \ll \nu\\
  \infty& \mbox{ else}
 \end{array}
 \right.\ee
 where $ \mu \ll \nu $  denotes that $\mu$ is absolutely continuous with respect to $\nu$. For any $f$ in $B(E)$, let $\langle \mu, f\rangle$ denote the integration of $f$ with respect to $\mu$. It  is known (\cite{DV75I}) that 
 \be\label{entropy-var}
 H(\mu|\nu)=\sup_{f\in B(E)}\{\langle \mu,f\rangle -\log\langle\nu, e^f\rangle\} =\sup_{f\in C(E)}\{\langle \mu,f\rangle -\log\langle\nu, e^f\rangle\}.  \ee
 
 For any $m\geq 1$, set
\[
E^m=\underbrace{E\times\cdots \times E}_m \]
and 
\[
\triangle_m =\{{\bf p}_m=(p_1,\ldots, p_m)\in E^m:   \sum_{i=1}^m p_i\leq 1\}.
\]
Define

\[
E^\infty=\underbrace{E\times\cdots \times E}_\infty\]

\[
\triangle_{\infty}=\{{\bf p}=(p_1, p_2, \ldots)\in E^{\infty}: \sum_{i=1}^{\infty}p_i\leq 1\}.
\]

Both $\triangle_m$ and $\triangle_{\infty}$ are equipped with the respective subspace topologies of $E^m$ and $E^{\infty}$.  We will use the following metrics that generate these topologies. 

\beq
d_m({\bf u}_m,{\bf v}_m)&=&\max\{|u_i-v_i|: i=1, \ldots, m\},\ {\bf u}_m, {\bf v}_m \in E^m\\
d({\bf u}, {\bf v})&=&\sum_{i=1}^{\infty}\frac{|u_i-v_i|\wedge 1}{2^i},\ \ {\bf u}, {\bf v} \in E^{\infty}.
\eeq

    
\section{Asymptotic Results}
\setcounter{equation}{0}

In this section, we establish the law of large numbers and the large deviation principles.

\subsection{Law of Large Numbers}

\begin{theorem}\label{lln}

Let ${\bf 0}= (0,0,\ldots)$ denote the origin in $\triangle_{\infty}$. As $\alpha$ and $\beta$ tend to infinity,  $\mathbb{Z}_{\alpha,\beta}$ and  $\Xi_{\alpha,\beta,\nu_0}$ converge in probability to ${\bf 0}$ and $\nu_0$ respectively.

\end{theorem}

\proof  
For any ${\bf p}, {\bf q}$ in $\triangle_{\infty}$, we have 

\[
d({\bf p},{\bf q})=\sum_{i=1}^{\infty}\frac{|p_i-q_i|}{2^i}.\]

By direct calculation, 
\beq
\mathbb{E}[Z_i]&=& \mathbb{E}[\mathbb{E}[Z_i |V_1, \ldots,V_i]] \\
&=&\mathbb{E}\bigg[\bigg(\prod_{k=1}^{i-1}\mathbb{E}[(1-W_k) |V_1, \ldots,V_i] \bigg) \mathbb{E}[W_i |V_1, \ldots,V_i] \bigg]\\
&=&\mathbb{E}\bigg[(1-V_1)\cdot\frac{1-V_1-V_2}{1-V_1}\cdots \frac{1-\sum_{k=1}^{i-1}V_k}{1-\sum_{k=1}^{i-2}V_k} \cdot\frac{V_i}{1-\sum_{k=1}^{i-1}V_k}\bigg]\\
&=& \mathbb{E}[V_i]=\bigg(\frac{\alpha}{1+\alpha}\bigg)^{i-1}\frac{1}{1+\alpha}\\
&=&\mathcal{O}(\frac{1}{\alpha})
\eeq

For any $\delta >0$, let $n_\delta$ be an integer such that $2^{-n_{\delta}}<\delta$. Then we have that for any $\ep >0$,
\beq
P\{d(\mathbb{Z}_{\alpha,\beta}, {\bf 0})\geq \ep\}&\leq & \ep^{-1}\left[\sum_{i=1}^{n_{\delta}}\mathbb{E}[Z_i]+\delta \right]\\
&=&\mathcal{O}(\frac{1}{\alpha})+\ep^{-1}\delta
\eeq
which converges to zero by taking the limit of $\alpha$ going to infinity followed by $\delta$ going to zero.

Next we turn to the limit of $\Xi_{\alpha,\beta,\nu_0}$. For each $f$ in $C(E)$, we have 
\[
\langle \Xi_{\alpha,\beta,\nu_0}, f\rangle=\sum_{i=1}^{\infty}Z_i f(\xi_i),\]

and

\beq
\mathbb{E}[\langle \Xi_{\alpha,\beta,\nu_0}, f\rangle]&=& \langle \nu_0, f\rangle\\
\vspace{3mm}
\mathbb{E}[\langle \Xi_{\alpha,\beta,\nu_0}, f\rangle^2]&=& \mathbb{E}[\sum_{i=1}^\infty Z_i^2] \langle \nu_0, f^2\rangle+  \mathbb{E}[\sum_{i\neq j}^\infty Z_i Z_j] \langle \nu_0, f\rangle^2\\
&=&   \mathbb{E}[\sum_{i=1}^\infty Z_i^2] [\langle \nu_0, f^2\rangle-\langle \nu_0, f\rangle^2] + \langle \nu_0, f\rangle^2.\eeq 

Similarly, by exploring the conditional beta structure, we have that for each $i\geq 1$, 

\beq
\mathbb{E}[Z^2_i]&=& \mathbb{E}[\mathbb{E}[Z^2_i |V_1, \ldots,V_i]] \\ 
&=&\mathbb{E}\bigg[\bigg(\prod_{k=1}^{i-1}\mathbb{E}[(1-W_k)^2 |V_1, \ldots,V_i] \bigg) \mathbb{E}[W^2_i |V_1, \ldots,V_i] \bigg]\\
&=&\mathbb{E}\bigg[\frac{[\beta(1-V_1)+1]\beta(1-V_1)}{\beta(1+\beta)}\cdot\frac{[\beta(1-V_1-V_2)+1]\beta(1-V_1-V_2)}{[\beta(1-V_1)+1]\beta(1-V_1)}\\
&&\hspace{-8mm} \cdots \frac{[\beta(1-\sum_{k=1}^{i-1}V_k)+1]\beta(1-\sum_{k=1}^{i-1}V_k)}{[\beta(1-\sum_{k=1}^{i-2}V_k)+1] \beta(1-\sum_{k=1}^{i-2}V_k)} 
\cdot\frac{(\beta V_i+1)\beta V_i}{[\beta(1-\sum_{k=1}^{i-1}V_k)+1]\beta(1-\sum_{k=1}^{i-1}V_k)}\bigg]\\
&=& \frac{1}{1+\beta}\mathbb{E}[V_i]+ \frac{\beta}{1+\beta}\mathbb{E}[V_i^2]\\
&=&\frac{1}{1+\beta}\bigg(\frac{\alpha}{1+\alpha}\bigg)^{i-1}\frac{1}{1+\alpha}+\frac{\beta}{1+\beta}\bigg(\frac{\alpha}{\alpha+2}\bigg)^{i-1}\frac{2}{(\alpha+2)(\alpha+1)}\bigg]. 
\eeq

 It follows that
 \[
 \mbox{Var}[\langle \Xi_{\alpha,\beta,\nu_0}, f\rangle]= \bigg[\frac{1}{1+\beta}+\frac{\beta}{1+\beta}\frac{1}{\alpha+1}\bigg]  [\langle \nu_0, f^2\rangle-\langle \nu_0, f\rangle^2]\]
 which converges to zero as $\alpha$ and $\beta$ tend to infinity.  Let $n_\delta$ be defined as above. Then

 \beq
 P\{\rho(\Xi_{\alpha,\beta,\nu_0}, \nu_0)>2\delta\}&\leq&  \sum_{i=1}^{n_{\delta}}P\{|\langle \Xi_{\alpha,\beta,\nu_0}-\nu_0, f_i\rangle|> \frac{2^{n_{\delta}}\delta}{n_{\delta}}\}  \\
 &\leq& \frac{n^2_\delta}{2^{n_{\delta}+1}\delta^2}\sum_{i=1}^{n_{\delta}} \mbox{Var}[\langle \Xi_{\alpha,\beta,\nu_0}, f_i\rangle]\ra 0, \  \alpha\ra \infty, \beta \ra \infty \eeq
 which leads to the law of large numbers for $\Xi_{\alpha, \beta, \nu_0}$.
 
 \hfill $\Box$


\subsection{Large Deviations}

The focus of  this subsection will be on the large deviations for $Q_{\alpha,\beta}$ and $\Pi_{\alpha,\beta, \nu_0}$. Due to the different topological structures, we prove the results separately by exploring the corresponding local structures.

\begin{theorem}\label{ldp-mass}

Assume that \rf{assump} holds. Then the family $\{Q_{\alpha,\beta}: \alpha>0, \beta>0\}$ satisfies a large deviation principle on space $\triangle_{\infty}$ with speed $\gamma =\max\{\alpha, \beta\}$ and good rate function
\be\label{mass-ratefunction}
I({\bf z})=\sup_{m\geq 1}I_m(z_1,\cdots, z_m)
 \ee
where
\beqn
&&I_m(z_1,\cdots, z_m)= \inf\bigg\{\sum_{i=1}^m \bigg(a \log \frac{1}{1-u_i}+b \prod_{j=1}^{i-1}(1-u_j)h(u_i,w_i)\bigg): \label{mass-ratefunction-finite}\\
&&\hspace{0.5cm} u_i,w_i \in E, u_i<1, i\geq 1, 
 (w_1, \cdots,  (1-w_1)\cdots(1-w_{m-1})w_m)=(z_1,\cdots, z_m)\bigg\}\nn
\eeqn 
and 

\beq
h(u, w) &=&  u\log \frac{u}{w}+(1-u)\log\frac{1-u}{1-w},\\
 \prod_{j=1}^{i-1}(1-u_j)&=& 1\ \ \mbox{for}\ i=1.
 \eeq

The coefficients $a$ and $b$ are given by
\[
(a,b)=\left\{\begin{array}{ll}
(c,1) & c<1\\
(1,c^{-1}) &  c>1\\
(1,1) &  c=1\\
\end{array}
\right.
\]
\end{theorem}

\proof Since the space $\triangle_{\infty}$ can be identified as the projective limit of $\triangle_m, m\geq 1$, by the Dawson-G\"artner theorem (\cite{DeZe98}), it suffices  to show that for each $m\geq 1$ the law of $(Z_1, \ldots, Z_m)$ satisfies a large deviation principle with speed $\ga$ and good rate function $I_m(\cdot)$. To do this we start with the large deviations for $({\bf U}_m,{\bf W}_m)=(U_1, \cdots, U_m, W_1,\cdots, W_m)$ and then apply the contraction principle.  Since the state space of $({\bf U}_m,{\bf W}_m)$ is compact, it follows from Theorem P in \cite{Puk91} that we only need to show the existence, and to obtain the expression of the limit

\beqn
&&\lim_{\delta \ra 0}\liminf_{\ga \ra \infty}\frac{1}{\ga}\log P\{({\bf U}_m, {\bf W}_m)\in B({\bf u}_m, {\bf w}_m; \delta)\}\label{local-ldp}\\
&&\ \ \ \ \ \ \ \ =\lim_{\delta \ra 0}\limsup_{\ga \ra \infty}\frac{1}{\ga}\log P\{({\bf U}_m, {\bf W}_m) \in \bar{B}({\bf u}_m, {\bf w}_m; \delta)\}\nn
\eeqn
for any ${\bf u}_m, {\bf w}_m$ in $E^m$, where 
\beq
B({\bf u}_m, {\bf w}_m; \delta)&=& \{{\bf x}_m, {\bf y}_m\in E^{m}: d_m({\bf x}_m, {\bf u}_m)<\delta, d_m({\bf y}_m, {\bf w}_m)<\delta\}\\
\bar{B}({\bf u}_m, {\bf w}_m; \delta)&=& \{{\bf x}_m, {\bf y}_m\in E^{m}: d_m({\bf x}_m, {\bf u}_m)\leq\delta, d_m({\bf y}_m, {\bf w}_m)\leq\delta\}\eeq

Since the function $\log x-\log y $ is not continuous at the origin, we need to  divide the discussion into several cases.
\vspace{2mm}
 
{\bf Case 1:} $u_i =1$ for some $i=1,\ldots,m$.  
\vspace{2mm}

By direct calculation, we have that 
\beq
&&\lim_{\delta \ra 0}\limsup_{\ga \ra \infty}\frac{1}{\ga}\log P\{({\bf U}_m, {\bf W}_m) \in \bar{B}({\bf u}_m, {\bf w}_m; \delta)\}\\
&&\ \ \ \ \ \ \ \  \leq \lim_{\delta \ra \infty}\limsup_{\ga \ra \infty}\frac{1}{\ga}\log P\{|U_i-u_i|\leq \delta\}\\
&&\ \ \ \ \ \ \ \  \leq \lim_{\delta \ra \infty}\limsup_{\ga \ra \infty}\frac{1}{\ga}\log P\{1-\delta\leq U_i \leq 1\}\\
&& \ \ \ \  \ \ \ \ \ \ \ \  =\lim_{\delta \ra 0}\limsup_{\ga \ra \infty}\frac{1}{\ga}\log \delta^{\alpha}=-\infty\eeq
which implies that \rf{local-ldp} holds with limit $-\infty$.
\vspace{2mm}

{\bf Case 2:}  $u_k < 1$ for all $1\leq k\leq m$ and $u_j=w_j= 0$ for some $j$.  
\vspace{2mm}

 For any $ i =1, \ldots, m$ and ${\bf x}_m$ in $E^m$, set $v_i=(1-x_1)\cdots(1-x_{i-1})x_i$ with $x_0=0$ and 
\[
 g_i(x_i, y_i; x_1,\ldots, x_{i-1})=(1-x_i)^{\alpha-1}\frac{ \Gamma(\beta(1-\sum_{k=1}^{i-1}v_k))}{\Gamma(\beta v_i)\Gamma(\beta (1-\sum_{k=1}^i v_k))}y^{\beta v_i-1 }_i(1-y_i)^{\beta(1-\sum_{k=1}^i v_k)-1}.\]

It follows from the definition that 
\beqn
&&P\{\bar{B}({\bf u}_m, {\bf w}_m; \delta)\}\nn \\
&&\ \ \ \ = \alpha^m \idotsint\limits_{\bar{B}({\bf u}_m, {\bf w}_m; \delta)}\prod_{i=1}^m  g_i(x_i,y_i;x_1, \ldots,x_{i-1})d\,x_i d\,y_i \label{local-e}\\
&&\ \ \ \   = \alpha^m \idotsint\limits_{\bar{B}({\bf u}_m, {\bf w}_m; \delta)}\exp\{-\ga\sum_{i=1}^m  \frac{1}{\ga}\log g^{-1}_i(x_i,y_i;x_1, \ldots,x_{i-1})\}d\,x_1\cdots d\,x_m d\,y_1\cdots d\,y_m.\nn
\eeqn

It is clear from the definition that
\be\label{zero-upper}
\int_{0}^{\delta}\int_{0}^{\delta} g_j(x_j,y_j; x_1, \ldots, x_{j-1})d\,x_j d\,y_j\leq 1\ee
uniformly for all $x_1, \ldots, x_{j-1}$.
  
By Stirling's formula we have that  for $\beta\geq 2$
\beqn
&&\int_0^{\delta}\int_0^{\delta} g_j(x_j,y_j; x_1, \ldots, x_{j-1})d\,x_j d\,y_j\nn\\
&&\ \ \ \ \geq  (1-\delta)^{\alpha+\beta-2}\int_0^{\delta}\int_0^{\delta}\frac{\Gamma(\beta(1-\sum_{k=1}^{j-1}v_k))}{\Gamma(\beta v_j)\Gamma(\beta (1-\sum_{k=1}^j v_k))}y_j^{\beta v_j-1 }d\,x_j  d\,y_j\nn\\
&&\label{zero-lower}\\
&&\ \ \ \  = (1-\delta)^{\alpha+\beta-2}\int_0^{\delta}  \frac{\Gamma(\beta(1-\sum_{k=1}^{j-1}v_k))}{\Gamma(\beta v_j+1)\Gamma(\beta (1-\sum_{k=1}^j v_k))}\delta^{\beta v_j}d\,x_j\nn\\
&& \ \ \ \ \geq c_0 (1-\delta)^{\beta-1}\delta^{\beta\delta+1}\sqrt{\frac{1-\delta}{\beta \delta +1}}e^{\beta r_{\delta}} \nn
\eeqn
uniformly for all $x_1, \ldots, x_{j-1}$, where $c_0$ denotes a generic positive constant and 
\[
r_{\delta} =\inf\{v\log \frac{v}{v+1/2}: 0\leq v\leq \delta\} \ra 0 \ \mbox{as}\  \delta \ra 0.
\]

Putting together \rf{zero-upper} and \rf{zero-lower} we conclude that  for $u_j=w_j=0$  the integration of the density $
g_j(x_j,y_j; x_1, \ldots, x_{j-1})$ makes zero contribution to the limits
\be\label{ul-bound}
\lim_{\delta\ra 0}\liminf_{\ga \ra \infty  }\frac{1}{\ga}\log P\{B({\bf u}_m, {\bf w}_m; \delta)\}\ 
\mbox{and}\  \lim_{\delta\ra 0}\limsup_{\ga \ra \infty  }\frac{1}{\ga}\log P\{\bar{B}({\bf u}_m, {\bf w}_m; \delta)\}.\ee

\vspace{2mm}

{\bf Case 3:} For all $1\leq k \leq m$, $u_k<1$, and $u_k$ and $w_k$ are not equal to zero at the same time.
\vspace{2mm}

We choose $\delta$ small enough so that $u_k +\delta <1$ for all $k$, and $ w_r-\delta>0, u_l-\delta>0$ for $u_l>0, w_r>0$.  For any $1\leq i \leq m$, we have that

\[
 -\frac{1}{\ga}\log g_i(x_i,y_i; x_1, \ldots, x_{i-1})= L_{i, \alpha,\beta}({\bf x}_m, {\bf y}_m)+\frac{\beta v_i+1}{\ga}\log (1+1/\beta v_i)+o(1/\ga).
\]
where 

\beq
L_{i, \alpha,\beta}({\bf x}_m, {\bf y}_m)&=&\frac{\alpha-1}{\ga}\log \frac{1}{1-x_i}  + \frac{1}{\ga}\log[y_i(1-y_i)] \\
&& \hspace{1cm} + \frac{\beta}{\ga}(1-\sum_{k=1}^{i-1}v_k)\bigg[x_i\log\frac{x_i}{y_i}+(1-x_i)\log\frac{1-x_i}{1-y_i}\bigg] \eeq

The term $\frac{\beta v_i+1}{\ga}\log (1+1/\beta v_i)$ is clearly non-negative and finite for ${\bf x}_m$ satisfying $d_m({\bf x}_m, {\bf u}_m)< \delta$.  On the other hand, on the set  $\{{\bf x}_m: d_m({\bf x}_m, {\bf u}_m)< \delta, u_i+\delta/2\leq x_i\leq u_i+\delta\}$ we have 
 \[
 \frac{\beta v_i+1}{\ga}\log (1+1/\beta v_i)\leq \frac{\beta v_i+1}{\ga}\log (1+2/\beta \delta),\]
which converges to zero as $\ga$ tends to infinity.  Set 
\[
C({\bf u}_m, {\bf w}_m; \delta)=\{({\bf x}_m, {\bf y}_m) \in \bar{B}({\bf u}_m, {\bf w}_m; \delta): u_i+\delta/2  \leq x_i \leq u_i+\delta\}.\]

Controlling  $P\{\bar{B}({\bf u}_m, {\bf w}_m; \delta)\}$ from below by  $P\{C({\bf u}_m, {\bf w}_m; \delta)\}$, it follows that  the limit  \rf{ul-bound}
 receives zero contribution from the term  $\frac{\beta v_i+1}{\ga}\log (1+1/\beta v_i)$.  The function $L_{i, \alpha,\beta}({\bf x}_m, {\bf y}_m)$ converges to infinity if $u_i >w_i=0$ or $u_i<w_i=1$. Otherwise it is continuous at $({\bf u}_m, {\bf w}_m)$ with a finite value. Thus the equality \rf{local-ldp} holds with limit
 \[
 \sum_{i=1}^m \bigg[a \log \frac{1}{1-u_i}+b \prod_{j=1}^{i-1}(1-u_j)h(u_i,w_i)\bigg] \]
and the theorem follows.

 \hfill $\Box$ 
 
 \begin{remark} Under the assumption \rf{assump}, $\alpha$ and $\beta$ converge to infinity at the same magnitude. The representation in \rf{mass-ratefunction-finite} seems to indicate that the impact of level one becomes stronger when $c$ increases. In particular, one would expect that if we let $c$ going to infinity, then the rate function $I(\boldsymbol{z})$ will converge to  $-\log (1-\sum_{i=1}^\infty z_i)$, the large deviation rate function for the mass  $\mathbb{V}_{\alpha}$ of the Dirichlet process {\rm (}Theorem~2.2 in {\cite{Feng07}}{\rm)}.  But this is not true. In fact, by choosing $u_i=0$ for all $i$, we obtain 
 \[
 I({\bf z})\leq - b \log (1-\sum_{i=1}^\infty z_i). \]   
 
 Choosing $u_i=w_i$ for all $i$ in \rf{mass-ratefunction-finite}, we obtain
\[
I({\bf z})\leq -a \log(1-\sum_{i=1}^{\infty}z_i).
\]

If we choose $u_1=w_1=z_1, u_i=0, i>1$, then it follows  that 

\[
I({\bf z})\leq -[a \log(1-w_1)+ b(1-w_1)\sum_{k=2}^{\infty}\log (1-w_k)].
\]

Putting all these together,  we obtain that for $\sum_{i=1}^\infty z_i<1$
\[
I({\bf z})< -\log(1-\sum_{i=1}^{\infty}z_i).
\]

Thus $I({\bf z})$ is in general strictly less than the large deviation rate function for  $\mathbb{V}_{\alpha}$. The gap is large for $c$ near zero or infinity. This gap represents the impact of the hierarchical structure. It is easier for HDP to make large deviations than the Dirichlet process from the limit.
  \end{remark}

Next we turn to the large deviations for $\Pi_{\alpha,\beta, \nu_0}$. 

\begin{theorem}\label{ldp-measure}
 The family $\{\Pi_{\alpha,\beta,\nu_0}: \alpha>0, \beta>0\}$ satisfies a large deviation principle on space $M_1(E)$ with speed $\ga$ and good rate function
\[
J(\mu)=\left\{\begin{array}{ll} 
\inf_ {\nu \in M_1(E),  \mathrm{supp}(\nu) \subset \mathrm{supp}(\nu_0)}\{a H(\nu_0|\nu)\}+b H(\nu|\mu)\}& \mathrm{supp}(\mu) \subset \mathrm{supp}(\nu_0)\\
+\infty & \mbox{else}
\end{array}
\right.
\]
where 
$a, b$ are the same as in Theorem~\ref{ldp-mass}, and $\mathrm{supp}(\cdot)$ denotes the topological support of element in $M_1(E)$.\end{theorem}

\proof  We prove the theorem in the case $\mathrm{supp}(\nu_0)=E$. The proof for general cases requires only minor adjustment.

Since  $\Xi_{\alpha,\beta, \nu_0}$ is the image of $(\Xi_{\alpha,\nu_0}, \Xi_{\alpha,\beta,\nu_0})$ under the continuous projection, the theorem follows from the contraction principle and  the large deviation result for  $(\Xi_{\alpha,\nu_0}, \Xi_{\alpha,\beta,\nu_0})$. Applying Theorem P in \cite{Puk91} again the latter holds if we can show that for any $\mu, \nu$ in $M_1(E)$

\beqn
&&\lim_{\delta \ra 0}\liminf_{\ga \ra \infty}\frac{1}{\ga}\log P\{(\Xi_{\alpha,\nu_0}, \Xi_{\alpha,\beta,\nu_0})\in B(\nu, \mu; \delta)\}\label{local-measure-ldp}\\
&&\ \ \ \ \ \ \ \ =\lim_{\delta \ra 0}\limsup_{\ga \ra \infty}\frac{1}{\ga}\log P\{(\Xi_{\alpha,\nu_0}, \Xi_{\alpha,\beta,\nu_0}) \in \bar{B}(\nu, \mu; \delta)\}\nn
\eeqn
where
\beq
 B(\nu, \mu; \delta)&=&\{(\tau,\varsigma)\in M_1(E)\times M_1(E): \rho(\tau, \nu)<\delta, \rho(\varsigma, \mu)<\delta\}\\
 \bar{B}(\nu, \mu; \delta)&=&\{(\tau,\varsigma)\in M_1(E)\times M_1(E): \rho(\tau, \nu)<\delta, \rho(\varsigma, \mu)\leq\delta\}. \eeq

Fix $\mu,\nu$ in $M_1(E)$ and  Set
\[
\mathcal{P}_{\nu,\mu}=\{\pi=(t_1,\ldots ,t_{m-1}): m\geq 2, 0<t_1<\ldots<t_{m-1}<1, \mu(\{t_i\})=\nu(\{t_i\})=0\ \mbox{for all}\ i\}.
\]
Each $\pi$ in $\mathcal{P}_{\nu,\mu}$ corresponds to the partition $[0, t_1), \ldots, [t_{m-1},1]$ of $E$ and the total number of intervals in the partition will be denoted by $|\pi|$. For a given partition $\pi$ with $|\pi|=m$, we write
\[
\pi(\mu)=(\mu([0,t_1)), \ldots, \mu([t_{m-1},1])) \in \triangle_{m}.  
\]

It follows from the variational formula \rf{entropy-var} that
\[
aH(\nu_0|\nu)+b H(\nu|\mu)=\sup\{a H(\pi(\nu_0)|\pi(\nu))+b H(\pi(\nu)|\pi(\mu)): \pi \in \mathcal{P}_{\nu,\mu}\}\]

Since the support of $\nu_0$ is $E$, it follows that  
\beq
J(\mu)&=&\inf_ {\nu \in M_1(E)}\{aH(\nu_0|\nu)+b H(\nu|\mu)\}\\
&=&   \inf_ {\nu \in M_1(E)}\sup_{\pi \in \mathcal{P}_{\nu,\mu}}\{a H(\pi(\nu_0)|\pi(\nu))+b H(\pi(\nu)|\pi(\mu))\}.\eeq

For each $\pi$ in $\mathcal{P}_{\nu,\mu}$ with $|\pi|=m$ and $\delta >0$, define 
\beq
 B_\pi(\nu, \mu; \delta)&=&\{(\tau,\varsigma)\in M_1(E)\times M_1(E): d_{m}(\pi(\tau), \pi(\nu))<\delta, d_{m}(\pi(\varsigma), \pi(\mu))<\delta\}\\
 \bar{B}_\pi(\nu, \mu; \delta)&=&\{(\tau, \varsigma)\in M_1(E)\times M_1(E): d_{m}(\pi(\tau), \pi(\nu))\leq\delta, d_m(\pi(\varsigma), \pi(\mu))\leq\delta\}. \eeq

Since the function $(\pi(\tau), \pi(\varsigma))$ is continuous at $(\nu,\mu)$, it follows that for any $\delta_1>0$ there exists $\delta>$ such that
\[
 \bar{B}(\nu, \mu; \delta)\subset  \bar{B}_\pi(\nu, \mu; \delta_1).\]
 
 On the other hand, for the given $\delta$ there exists $M\geq 1$ such that
 \[
 \{(\tau, \varsigma) \in M_1(E)\times M_1(E):\sup_{1\leq i \leq M}\{|\langle\tau-\nu, f_i \rangle|\vee |\langle\varsigma-\mu, f_i \rangle|\}<\frac{\delta}{2M}\}
 \]
 is a subset of $B(\nu, \mu; \delta)$. Since $f_i$ is continuous for all $i$, it follows that
 there exists a partition $\tilde{\pi}$ and $\delta_2>0$ such that
 \[
 B_{\tilde{\pi}}(\nu, \mu; \delta_2)\subset B(\nu, \mu; \delta).\]
 
 Putting all these together it follows that \rf{local-measure-ldp} will hold if for each $\pi$ in $\mathcal{P}_{\nu,\mu}$
 \beqn
&&\lim_{\delta \ra 0}\liminf_{\ga \ra \infty}\frac{1}{\ga}\log P\{(\pi(\Xi_{\alpha,\nu_0}), \pi(\Xi_{\alpha,\beta,\nu_0}))\in B_{\pi}(\nu, \mu; \delta)\}\nn
\\
&&\ \ \ \ \ \ \ \ =\lim_{\delta \ra 0}\limsup_{\ga \ra \infty}\frac{1}{\ga}\log P\{(\pi(\Xi_{\alpha,\nu_0}), \pi(\Xi_{\alpha,\beta,\nu_0})) \in \bar{B}_{\pi}(\nu, \mu; \delta)\}   \label{local-finite-measure-ldp} \\
&& \ \ \ \ \ \ \ \ = - [aH(\pi(\nu_0)|\pi(\nu))+b H(\pi(\nu)|\pi(\mu))].\nn
\eeqn 

For a given partition $\pi=(t_1, \ldots,t_{m-1})$ in $\mathcal{P}_{\nu, \mu}$, denote $\pi(\mu)$ and $\pi(\nu)$ by ${\bf u}_m$ and ${\bf v}_m$ respectively, and let ${\bf r}_m=(r_1, \ldots,r_m)=\pi(\nu_0)$. Since $\mathrm{supp}(\nu_0)=E$, it follows that $r_i>0$ for all $i$.  The joint density function of $ (\pi(\Xi_{\alpha,\nu_0}), \pi(\Xi_{\alpha,\beta,\nu_0}))$ is

\beq
F({\bf q}_m, {\bf p}_m)&=&\frac{\Gamma(\alpha)}{\Gamma(\alpha r_1)\cdots\Gamma(\alpha r_m)}\frac{\Gamma(\beta)}{\Gamma(\beta q_1)\cdots\Gamma(\beta q_m)}\\
&& \hspace{1cm} \times \prod_{i=1}^m q_i^{\alpha r_i-1}p_i^{\beta q_i -1}, {\bf q}_m, {\bf p}_m \in \triangle_m.
\eeq

If $m=2$ and  $u_1=0$ (the case $u_1=1$ is similar),  then by Stirling's formula we have 
\beq
&&P\{(\pi(\Xi_{\alpha,\nu_0}), \pi(\Xi_{\alpha,\beta,\nu_0})) \in \bar{B}_{\pi}(\nu, \mu; \delta)\}\\
&&= \idotsint\limits_{\bar{B}_{\pi}(\nu, \mu; \delta)} F({\bf q}_2, {\bf p}_2) d\,q_1 d\,p_1 \\
&&= \frac{\Gamma(\alpha)}{\Gamma(\alpha r_1)\Gamma(\alpha r_2)}
 \int_{(v_1-\delta)\vee 0}^{(v_1+\delta)\wedge 1}q_1^{\alpha r_1-1}(1-q_1)^{\alpha r_2-1 }d\,q_1\\
 && \hspace{1cm}\times \int_0^\delta   \frac{\Gamma(\beta)}{\Gamma(\beta q_1)\Gamma(\beta q_2)}  p_1^{\beta q_1-1}(1-p_1)^{\beta q_2-1}d\,p_1\\
 &&\leq \frac{\Gamma(\alpha)}{\Gamma(\alpha r_1)\Gamma(\alpha r_2)}
 \int_{(v_1-\delta)\vee 0}^{(v_1+\delta)\wedge 1} \frac{\beta \delta^{\beta q_1}(1-\delta)^{\beta q_2}\Gamma(\beta)}{(1-\delta)^{\beta q_2+1}\Gamma(\beta q_1+1)\Gamma(\beta q_2+1)} q_1^{\alpha r_1-1}(1-q_1)^{\alpha r_2}d\,q_1\\
 &&\leq \int_{(v_1-\delta)\vee 0}^{(v_1+\delta)\wedge 1}\exp\{-\ga[\alpha/\ga H({\bf r}_2|{\bf q}_2)+\beta/\ga H({\bf q}_2|\boldsymbol{\delta}_2)+o(1/\ga)]\}d\,q_1
 \eeq
where $\boldsymbol{\delta}_2=(\delta,1-\delta)$.  If $ 0<v_1<1$, then $ H({\bf q}_2|\boldsymbol{\delta}_2)$ converges to infinity as $\delta$ tends to zero. If $v_1=0 \ \mbox{or}\  1$, then $ H({\bf r}_2|{\bf q}_2)$ converges to infinity as $\delta$ tends to zero.  Thus

 \be\label{upper-2-zero}
 \lim_{\delta \ra 0}\limsup_{\ga \ra \infty}\frac{1}{\ga}\log P\{(\pi(\Xi_{\alpha,\nu_0}), \pi(\Xi_{\alpha,\beta,\nu_0})) \in \bar{B}_{\pi}(\nu, \mu; \delta)\} =-\infty \ee 
and \rf{local-finite-measure-ldp} holds.

If $m>2$ and there exists $1\leq i\leq m$ such that $u_i=0$, then by the partition property of the Dirichlet process we can separately amalgamate all zero terms, and all non-zero terms to get  a partition with $m=2$ and use the above argument to show that  \rf{local-finite-measure-ldp} holds.

It remains to prove the result for the case $m\geq 2, u_i>0$ for all $i$. By Stirling's formula, we have 

\be\label{stir-1}
\log \frac{\Gamma(\alpha)}{\Gamma(\alpha r_1)\cdots\Gamma(\alpha r_m)}=\log \prod_{i=1}^m r_i^{-\alpha r_i} +o(\alpha)
\ee
and 
\beqn
\log\frac{\Gamma(\beta)}{\Gamma(\beta q_1)\cdots\Gamma(\beta q_m)}&=&\log \frac{\beta^m \prod_{i=1}^m q_i\Gamma(\beta)}{\Gamma(\beta q_1+1)\cdots\Gamma(\beta q_m+1)}\label{stir-2}\\
&=&\log \prod_{i=1}^m q_i^{-\beta q_i} + \sum_{i=1}^m (\beta q_i+1)\log \frac{\beta q_i}{\beta q_i+1} +o(\beta).\nn
\eeqn

Putting these together we obtain
\beq
\log F({\bf q}_m, {\bf p}_m)&=& -\ga[\alpha/\ga H({\bf r}_m|{\bf q}_m)+\beta/\ga H({\bf q}_m|{\bf p}_m) \\
&&- \beta/\ga \sum_{i=1}^m q_i\log \frac{ q_i}{q_i+1/\beta}+ \ga^{-1}\sum_{i=1}^m \log  p_i+o(1/\ga)].
\eeq

Since $u_i>0$ for all $i$, both the term $ \beta/\ga \sum_{i=1}^m q_i\log \frac{ q_i}{q_i+1/\beta}$ and  the term  $\ga^{-1}\sum_{i=1}^m \log  p_i$ converge to zero as $\ga$ goes to infinity followed by $\delta$ going to zero.  The result then follows from the fact that the term $\alpha/\ga H({\bf r}_m|{\bf q}_m)+\beta/\ga H({\bf q}_m|{\bf p}_m) $ converges to $ aH(\pi(\nu_0)|\pi(\nu))+b H(\pi(\nu)|\pi(\mu))$.

\hfill $\Box$ 

\begin{remark}\label{re-2}
It is known {\rm (}\cite{LS87}, \cite{DF01}{\rm)} that the large deviation rate function for the Dirichlet process $\Xi_{\alpha,\nu_0}$ is given by $H(\nu_0|\mu)$. 
Choosing $\nu=\mu\ \mbox{or}\ \nu_0$ we obtain that
\[
a H(\nu_0|\nu)+b H(\nu|\mu)= a H(\nu_0|\mu)\ \mbox{or} \ b H(\nu_0|\mu)\]
which implies that
\[
J(\mu) \leq \min\{a H(\nu_0|\mu), b H(\nu_0|\mu) \}.
\]

Thus the rate function for the HDP is less than the rate function for the Dirichlet process, which is consistent with the fact that  the number of clusters under HDP has a slower growth rate than under the Dirichlet process. It is not clear how to get a more explicit and simpler form for $J(\mu)$ even for concrete $\nu_0$.  
\end{remark}

\begin{remark}\label{re-3}
Generalization to HDP with more than two levels  will lead to a similar rate function involving the relative entropies of all levels. 
\end{remark}

\vspace{3mm}

\noindent \ack{ The author wishes to thank the referee for the careful review of the paper and many insightful suggestions.}

\bibliographystyle{amsplain}





\end{document}